# On estimates of trigonometrical integrals.

Jabbarov I. Sh.

### ABSTRACT

In this article a new upper bounds for the multiple trigonometrical integrals are found. The method of the work based on a new method of estimation for the areas of algebraic surfaces.

## 1. Introduction

As a multiple trigonometrical integral we call an integral

$$\int_{\Omega} G(\overline{x}) e^{2\pi i F(\overline{x})} d\overline{x}, \qquad (1)$$

where $\Omega$ denote some domain of $n$ dimensional space $R^n$; over the functions $G(x)$ and $F(x)$ one imposes definite conditions on boundedness or smoothness. Many researches (see [1-12, 30-33, 35-39]) are devoted to a question on estimates of trigonometrical integrals. The first results in this direction belong Van der Korput and E.Landau (see [33]). The result established in the work [4] where the authors have received a non improvable estimation for trigonometrical integrals has important applications. The multidimensional case also was investigated in the literature. Unlike to the one-dimensional case, estimating of multiple trigonometrical integrals of a view (1), where $\Omega$- some Jordan domain with a smooth boundary, and the functions $G(\overline{x}), F(\overline{x})$ are from a certain class of smoothness, much more difficult.

The scheme of receiving estimates of integrals of a view (1) is similar to the scheme of one-dimensional case. After some transformations the integral reduces to the view

$$\int_{a}^{b} V(u) e^{iu} du,$$

where $V(u)$ represents the surface integral depending on parameter $u$.

During more detailed consideration the difficulties of two types arise. The first consists in dissecting of the interval $[a,b]$ into the finite number of subintervals, in each of which $V(u)$ is monotonic, or in dissecting the domain $\Omega$ into finite number of such subdomains that in every of them the function $V(u)$ is monotonic. Other task consists in an estimation of some areas that allows us to estimate $\max|V(u)|$. These questions are solved in the present article. Below we prove general theorems on estimations of areas of some surfaces and apply them to the estimates of trigonometrical integrals for the class of polynomials standing on the exponent containing the polynomials a senior form of which contains all independent variables and which are not represented in the form of the sum of various polynomials in different varables. We consider here only the simplest integrals of a view

$$\int_\Omega e^{2\pi i F(\overline{x})} d\overline{x} . \qquad (2)$$

## 2. Formulation of results

Let $\Omega$ be a bounded closed domain of $n$-dimensional space $R^n$, $n \geq 1$, where $R$ denote a real line. Let's assume that in $\Omega$ an $r$-dimensional surface be given by means of a system of polynomial equations

$$f_j(\overline{x}) = 0, \; j = 1,...,n-r, 0 \leq r \leq n, \quad (3)$$

with a Jacoby matrix

$$J = J(\overline{x}) = \left\| \frac{\partial f_j}{\partial x_i} \right\|, \quad i = 1,...,n, \; j = 1,...,n-r ,$$

having, everywhere in $\Omega$, the maximal rank.

Let, further, $A_0 = A_0(\overline{x})$ be other functional matrix written down in a form

$$A_0 = A_0(\overline{x}) = \left\| f_{ij}(\overline{x}) \right\|, 1 \leq i \leq r, 1 \leq j \leq m$$

with polynomial entries. Having elements of columns of a matrix $A_0$ arranged in a line

$$f_{11}(\overline{x}),...,f_{r1}(\overline{x}), f_{12}(\overline{x}),...,f_{r2}(\overline{x}),...,f_{1m}(\overline{x}),...,f_{rm}(\overline{x})$$

let's take the transposed Jacoby matrix of this system of functions, designating it $A_1$:

$$A_1 = A_1(\overline{x}) = \begin{Vmatrix} \dfrac{\partial f_{11}}{\partial x_1} & \cdots & \dfrac{\partial f_{r1}}{\partial x_1} & \cdots & \dfrac{\partial f_{1m}}{\partial x_1} & \cdots & \dfrac{\partial f_{rm}}{\partial x_1} \\ \cdots & \cdots & \cdots & \cdots & \cdots & \cdots & \cdots \\ \dfrac{\partial f_{11}}{\partial x_n} & \cdots & \dfrac{\partial f_{r1}}{\partial x_n} & \cdots & \dfrac{\partial f_{1m}}{\partial x_n} & \cdots & \dfrac{\partial f_{rm}}{\partial x_n} \end{Vmatrix}.$$

Then, elements of columns of this matrix consequently as above, we arrange in a line and take the transposed Jacoby matrix $A_2 = A_2(\overline{x}) = A_1'(\overline{x})$ of the received system of functions, and we will continue this procedure while we have not received a matrix $A_k = A_{k-1}'(\overline{x})$ for a given $k \geq 1$. The last matrix defined thus consists of every possible partial derivatives of the same order $k$ of elements of the matrix $A_0 = A_0(\overline{x})$ and has the size $n \times n^{k-1}rm$. Let's assume that $A_j(\overline{x})$ has in $\Omega$ the maximal rank equal $n$. Let's designate $G_j(\overline{x})$ the product of the last (smallest) $r$ singular numbers of the matrix $A_j(\overline{x})$, $j = 0,...,k$. We put

$$E = E(H) = \left\{ \overline{x} \in \Omega \big| G_0(\overline{x}) \leq H \right\}, H > 0.$$

If $\varphi_{ik}(\overline{x})$ are entries of the matrix $A_j(\overline{x})$ we will accept the following designations

$$L_j(\overline{x}) = \left( \sum_{i,k} \left| \varphi_{ik}(\overline{x}) \right|^2 \right)^{1/2},$$

$$L = \max_j \max_{x \in \Omega} L_j(\overline{x}), \quad G_j = \min_{x \in \Omega} G_j(\overline{x}), j = 0,...,k.$$

Below we will prove two auxiliary theorems allowing to estimate areas of a special type of surfaces arising at studying of questions connected with trigonometrical integrals.

For the preparation of formulation and the proof of our statements we need to dissect the domain $\Omega$ into such parts that on each of them the system (3) would allow one sheeted and one valued solvability. These parts are defined by maximal minors for the modulus of the Jacoby matrix $J$ of this system. Let's dissect $\Omega$ into no more than $t = C_n^{n-r}$ subdomains $\Omega^{(v)}, v = 1,...,t$, intersecting each with another only by parts of their boundary, in each of which one of minors of the Jacoby matrix has the maximal absolute values among all minors. In each subdomain of a view $\Omega^{(v)}$ the $v$-th minor (in some their ordering) everywhere accepts maximal absolute values. As it is well known each subdomain $\Omega^{(v)}$ is a closed set and *can be represented as a union of finite number of one-connected closed subdomains, as a set of solutions in $\Omega$ of the system of polynomial inequalities.* Therefore, each subdomain $\Omega^{(v)}$ is represented in the form of a union $\Omega^{(v)} = \bigcup_{c \leq T_0} \Omega(v, c)$, where $\Omega(v, c)$ is a one connected subdomain.

Let's consider a subdomain $\Omega(v, c)$. In this subdomain the system allows, in some neighbourhood of an arbitrary solution of the system (3), unique solvability in regard to *the same variables.* Each subdomain $\Omega(v, c)$ can be dissected into no more than $T$ subdomains $\Delta_\mu, \mu = 1, ..., f, f \leq F$ in each of which the system (3) allows one sheeted solvability concerning $n - r$ variables. Let $\bar{\xi} = (\xi_1, ..., \xi_r)$ be a vector of free variables. Then, $x_i$ is possible to present as a function $x_i = x_i(\bar{\xi})$ of free variables. Let's designate $A_0(\bar{\xi})$ the matrix received from the matrix $A_0(\bar{x})$ replacing the variables $x_i$ by functions $x_i = x_i(\bar{\xi})$. In other words we consider the functional matrix $A_0$ as a matrix depending on $\bar{\xi}$. Let's designate $G_{(1)}$ the minimal value of Gram determinant of gradients of elements $A_0(\bar{\xi})$ (differentiation is taken with regard to $\bar{\xi}$), i.e. $G_{(1)} = \min \det\left(A_{1\bar{\xi}} \cdot {}^t A_{1\bar{\xi}}\right)$ (we will notice that the minimum undertaken over all $\bar{\xi}$, $v$ and $c$). Thus $A_{1\bar{\xi}}$ means the matrix of the size $r \times rm$ received from $A_0$ by differentiation with regard to $\bar{\xi}$, i.e $A_{1\bar{\xi}} = A_0'(\bar{\xi})$. Therefore, the matrix $A_1(\bar{x})$ being considered as a matrix of $\bar{\xi}$, differs from $A_{1\bar{\xi}}$. Further, for positive number $a$ we designate $h(a) = a + a^{-1}$. It is obvious that $a \leq h(a)$, $h(a^{-1}) = h(a)$, and $h(ab) \leq h(a)h(b)$, at $a, b > 0$.

**Theorem 1**. Let $\Pi_H$ be a part of a surface (3) included in $E$ (H) and $G_{(1)} > 0$. Then for the area $\mu(\Pi_H)$ we have the bound

$$\mu\left(\Pi_H\right) \leq FT_0 \cdot 2^{r+3} r^{3r} c_0^2 \binom{nr}{r}^{1/2} \binom{n}{r}^{3/2} H \cdot G_{(1)}^{-1} \cdot \tilde{\wp}^r$$

$$\tilde{\wp} = r^2 \log\left\{h\left(G_{(1)}\right) h(H) h(L)\right\}, c_0 = \pi^{-r/2} \Gamma\left(1 + r/2\right).$$

Similarly, we can, beginning from the matrix $A_{j-1}$, form a matrix $A_{j\bar{\xi}}$ assuming that $G_{(j)} > 0$ for all considered $j > 0$.

**Theorem 2**. *Let $k \geq 1$ and $G_{(k)} > 0$. Then under the conditions of the theorem 1 we have:*

$$\mu\left(\Pi_H\right) \leq c_0^2 \left(1 + 4n^r c''\right)^{k-1} 2^{r+3} \left(nr^2\right)^{2^r} F T_0 H^{1/k} G_{(k)}^{-1/k} \tilde{\wp}_k^r;$$

$$\tilde{\wp}_k = 3r^2 \log \tilde{H}; \tilde{H} = \max\left\{h(H), h\left(G_{(1)}\right), \dots, h\left(G_{(k)}\right), h(L)\right\},$$

*and $c''$- some constant, and numbers $G_{(1)}, \dots, G_{(k)}$ are defined by equalities*

$$G_{(1)} = H^{1/2} G_{(2)}^{1/2}, G_{(2)} = H^{1/3} G_{(3)}^{1/3}, \dots, G_{(k)} = H^{1/k} G_{(k)}^{(k-1)/k}.$$

**Theorem 3**. *Let $k \geq 1$ and $G_k > 0$. Then under the conditions of the theorem 2:*

$$\mu\left(\Pi_H\right) \leq c_0^2 \left(1 + 4n^r c''\right)^{k-1} 2^{r+3} \left(nr^2\right)^{2^r} F T_0 H^{1/k} G_k^{-1/k} \tilde{\wp}_k^r.$$

Let $F(\bar{x})$ be some polynomial. Let's consider the trigonometrical integral (2), in the domain $\Omega$ with a piecewise algebraic boundary. As a matrix $A_0$ we take the matrix

$${}^t\nabla F = \left\|\frac{\partial F}{\partial x_1}, \dots, \frac{\partial F}{\partial x_r}\right\|.$$

Let everywhere in $\Omega$

$$\|\nabla F\| = \sqrt{\left(\frac{\partial F}{\partial x_1}\right)^2 + \dots + \left(\frac{\partial F}{\partial x_r}\right)^2} \neq 0.$$

The matrix $A_1(\bar{x})$ looks like:

$$\left\|\begin{matrix} \dfrac{\partial^2 F}{\partial x_1^2} \cdots \dfrac{\partial^2 F}{\partial x_1 \partial x_r} \\ \dfrac{\partial^2 F}{\partial x_r \partial x_1} \cdots \dfrac{\partial^2 F}{\partial x_r^2} \end{matrix}\right\|,$$

and the matrix $A_{k-1}(\bar{x})$ is combined of all partial derivatives of order $k$ of the function $F(\bar{x})$.

**Theorem 4**. *1) If $k \geq r$ then there exist a positive constant $c_1 = c_1(r,k,degF)$ such that*

$$\left| \int_{\Omega} e^{2\pi i F(\overline{x})} d\overline{x} \right| \le c_1 K G^{\frac{r-k}{k(k-1)}} G_{k-1}^{-\frac{r}{k(k-1)}} \left\{ log(L_0 + L_0^{-1}) \right\}^{r-1},$$

where $G = \min_{x \in \Omega} \sqrt{\det(A_{k-1} \cdot {}^t A_{k-1})}$, $G_{k-1}$ was is defined above;

2) If $k < r$ then the estimation

$$a) \left| \int_{\Omega} e^{2\pi i F(\overline{x})} d\overline{x} \right| \le c_2 K \widetilde{H}^{\frac{r-k}{k-1}} G_{k-1}^{-\frac{1}{k-1}} \left\{ log(L_0 + L_0^{-1}) \right\}^{r+1},$$

takes plase, where $c_2 = c_2(r, k, F)$ is a constant,

$$\widetilde{H} = \max_{\overline{x} \in \Omega} \| \nabla F \|, \quad H_1 = \min_{x \in \Omega} \| \nabla F \|, L_0 = \max\left(L, L^{-1}, \widetilde{H}, H_1^{-1}, G_{k-1}, G_{k-1}^{-1}\right);$$

$$b) \left| \int_{\Omega} e^{2\pi i F(\overline{x})} d\overline{x} \right| \le C_3 K \Pi^{\frac{r-k}{k-1}} G_{k-1}^{-\frac{1}{r-1}} \left\{ log(L_0 + L_0^{-1}) \right\}^{r-1},$$

where $\Pi$ denote the maximal area of the surface $F(\overline{x}) = u$ (the maximum is taken over $u$).

### 3. Auxiliary statements.

We assume that the boundary of the domain $\Omega$ defined by finite number of equations of a view $H_i(\overline{x}) = a, a \in R$ (the parameter $a$ can accept different values for different $i$) with polynomial functions on the left sides, thus the domain $\Omega$ becomes covered by a union of a finite family of surfaces of a view $H_i(\overline{x}) = c$, where the parameter $c$ accepts continuously all values from some seqment. Without breaking a generality, it is possible to count that the number of such functions is 1 and $H_1(\overline{x}) = H(\overline{x})$.

**Lemma 1.** There exist such a dissection of the domain $\Omega$ into the union of no more than finite number of subdomains so that the surface integral $\varphi(u) = \int_{F(\overline{x}) = u} \frac{ds}{\| \nabla F \|}$, respectively, breaks into the sum of the surface integrals being monotonous functions of a variable $u$, moreover, the number of addends the last sum depends only on degree of a polynomial $F$.

*Proof.* Having given to the variable $u$ some increment $\Delta u$, we can write

$$\varphi(u + \Delta u) - \varphi(u) = \int_{F(\overline{x}) = u + \Delta u} \frac{ds}{\| \nabla F \|} - \int_{F(\overline{x}) = u} \frac{ds}{\| \nabla F \|}.$$

As the domain $\Omega$ is closed, the gradient of function and its partial derivatives of the second order are bounded. Let's take for the function $F$ Taylor decomposition in a neighbourhood of the point $\bar{x}$, lying on the surface $F = u$, in the gradient direction:

$$F(\bar{x} + \lambda \nabla F) - F(\bar{x}) = \lambda(\nabla F, \nabla F) + o(\lambda^2).$$

Let's pick up $\lambda$ so that the point $\bar{x} + \lambda \nabla F$ was on the surface $F = u + \Delta u$. Then, we should have:

$$\Delta u = \lambda \|\nabla F\|^2 + o(\lambda^2).$$

At rather small $\Delta u$ the second term on the right part is small in comparison with the first. Therefore, we have

$$\lambda = \frac{\Delta u}{\|\nabla F\|^2} + o(\Delta^2 u).$$

Let's find now an increment $\delta$ of the function $\|\nabla F\|^{-1}$ during transition from the point $\bar{x}$ to the point $\bar{x} + \lambda \nabla F$ in the gradient direction. We have:

$$\delta = \left\| F(\bar{\xi} + \lambda \nabla F) \right\|^{-1} - \left\| F(\bar{\xi}) \right\|^{-1} = \nabla\left( 1 \middle/ \sqrt{(\partial F/\partial x_1)^2 + \cdots + (\partial F/\partial x_n)^2} \right) \cdot \lambda \nabla F(1 + o(\lambda)),$$

where on the right part of this equality a scalar product of gradients of the specified functions stands. We have:

$$\frac{\partial \|\nabla F\|^{-1}}{\partial x_i} = -\|\nabla F\|^{-3}\left( \sum_{j=1}^{n} \frac{\partial F}{\partial x_j} \cdot \frac{\partial^2 F}{\partial x_j \partial x_i} \right).$$

Then,

$$\nabla\left( 1 \middle/ \sqrt{(\partial F/\partial x_1)^2 + \cdots + (\partial F/\partial x_n)^2} \right) \cdot \lambda \nabla F = -\lambda \|\nabla F\|^{-3} \sum_{i=1}^{n} \sum_{j=1}^{n} \frac{\partial^2 F}{\partial x_j \partial x_i} \cdot \frac{\partial F}{\partial x_j} \frac{\partial F}{\partial x_i}.$$

Therefore, substituting the found value of $\lambda$, we find:

$$\delta = -\left( 1/\|\nabla F\|^3 \right)\left( A_1(\nabla F), \lambda \nabla F \right)(1 + o(\lambda^2)) =$$

$$= -\Delta u \|\nabla F\|^{-3}\left( A_1 \overline{\nabla}, \overline{\nabla} \right) + o(\Delta^2 u); \ \overline{\nabla} = \nabla F / \|\nabla F\|,$$

here $A_1 = \left( \dfrac{\partial^2 F}{\partial x_i \partial x_j} \right)$ is the matrix defined above, $\left( A_1 \overline{\nabla}, \overline{\nabla} \right)$ means a quadratic form.

Under the conditions imposed on a gradient, as shown above, the domain $\Omega$ may be dissected into finite number of subdomains which pairwisely intersecting only by parts of the boundary, and where the equation $F = u$ allows one sheeted solvability with respect one and the same variables. Let's consider one of them where the mentioned equation is solved with respect, say, to $x_1$:

$$x_1 = \varphi(x_2,...,x_n); \ (x_2,...,x_n) \in \omega ,$$

and $\omega$ is an domain of changing for independent variables. Having taken the any point $\bar{\xi} \in \omega$, we will define the mapping $\varphi$ putting to each point $\bar{\xi} \in \omega$ in correspondance a point $\left( \varphi\left(\bar{\xi}\right), \bar{\xi} \right)$ on the surface $F = u$, and will consider tangential linear mapping

$$\Phi : \bar{\xi} \mapsto \varphi\left(\bar{\xi}\right) + \varphi'\left(\bar{\xi}\right) \cdot \Delta\bar{x}; \ \bar{\xi} = (\xi_1,...,\xi_n) \in \omega . \ \ (4)$$

The image of this mapping is a tangential linear variety (hyper plane) to the surface $F = u$ in the point $\left( \varphi\left(\bar{\xi}\right), \bar{\xi} \right)$. Let's notice that the point $\left( \Phi\left(\bar{\xi}\right), \bar{\xi} \right)$ of the tangential hyper plane will situated from the corresponding point $\left( \varphi\left(\bar{\xi}\right), \bar{\xi} \right)$ on the surface $F = u$ at a distance $o\left( \left| \Phi(\bar{x}) - \varphi(\bar{\xi}) \right| \right)$ which is of order $o(\Delta x)$. In each point $\bar{x}$ of the surface $F = u$ the gradient $\nabla F$ is orthogonal to the tangential hyper plane. Really,

$$\nabla F \cdot \Phi'(\xi)\Delta x = \left( \frac{\partial F}{\partial x_1},...,\frac{\partial F}{\partial x_n} \right) \cdot \left( \left( \frac{\partial F}{\partial x_1} \right)^{-1} \left[ -\frac{\partial F}{\partial x_2}\Delta x_2 - \cdots - \frac{\partial F}{\partial x_n}\Delta x_n \right], \Delta x_2, \cdots, \Delta x_n \right) = 0 .$$

When $\lambda$ is defined as above, the point $\bar{x} + \lambda\nabla F$ where $\bar{x} \in \Pi(u)$, belongs to the surface $\Pi(u + \Delta u)$; here by $\Pi(u)$ we designate the surface defined by the equation $F = u$ in a wider open domain $\Omega' \supset \Omega$. For any open domain $\Omega'$ the surface $\Pi(u + \Delta u) \cap \Omega$ entirely lies in $\Omega'$ for all enough small values of $|\Delta u|$. The mapping $\Psi : \Pi(u) \rightarrow \Pi(u + \Delta u)$ defined as $\Psi(\bar{x}) = \bar{x} + \lambda\nabla F$ is one-one mapping when $|\Delta u|$ is sufficiently small. Really,

$$\Psi(\bar{x}) = \bar{x} + \left( \frac{\Delta u}{\|\nabla F\|^2} + o(|\Delta u|) \right) \nabla F = \bar{x} + \Delta u \frac{\nabla F}{\|\nabla F\|^2} + o(|\Delta u|) ,$$

and at sufficiently small $|\Delta u|$ the Jacoby matrix of this mapping can be represented as a sum of identity matrix and a Jacoby matrix of the mapping

$$\bar{x} \mapsto \Delta u (1 + o(1)) \frac{\nabla F(\bar{x})}{\left\| \nabla F(\bar{x}) \right\|^2} \, .$$

Owing to the imposed constraints on $F(\bar{x})$ and on the domain $\Omega$, determinant of the Jacoby matrix $\Psi$ tends to 1 as $\Delta u \to 0$, i.e. this determinant will be distinct from zero everywhere in considered domain. So, $\Psi$ is a bijectiv mapping at sufficiently small $|\Delta u|$.

We put: $D(u) = \{ \bar{x} \in \Omega \mid F(\bar{x}) = u \}$. Then, the surface $D(u + \Delta u)$ tends to $D(u)$ as $\Delta u \to 0$. $\Psi(D(u))$ is a closed subset of $\Pi(u + \Delta u)$. Further, a prototype $D(u + \Delta u)$ of the same mapping we will designate $D'(u + \Delta u)$. Then, we have:

$$\varphi(u + \Delta u) - \varphi(u) = \int_{D'(u+\Delta u) \cap D(u)} \left( \frac{1}{\left\| \nabla F(\bar{x} + \lambda \nabla F) \right\|} - \frac{1}{\left\| \nabla F(\bar{x}) \right\|} \right) ds +$$

$$+ \int_{D(u+\Delta u) \setminus \Psi(D(u))} \frac{ds}{\left\| \nabla F \right\|} - \int_{D(u) \setminus D'(u+\Delta u)} \frac{ds}{\left\| \nabla F \right\|} \, . \qquad (4)$$

Substituting the value found above for an increment $\delta$, we find for the first surface integral the following expression:

$$- \Delta u (1 + o(1)) \int_{F(\bar{x})=u} \left\| \nabla F \right\|^{-3} \left( A_1 \nabla, \nabla \right) ds$$

Let's consider two remained surface integrals. They will be transformed equally. The first integral is taken over the surface $D(u + \Delta u) \setminus \Psi(D(u))$ which is included between the boundaries $D(u + \Delta u)$ and $\Psi(D(u))$. It is clear that as $\Delta u \to 0$ then this piece is narrowing, being pulled off along $n-2$-dimensional surface of intersection $D(u + \Delta u) \cap \partial \Omega$, which in turn, tends to the limiting position $D(u) \cap \partial \Omega$ (it may be empty).

Rerecall the construction of the surface integral given in [17, p. 261]. Let $\omega'$ be an $n-1$-dimensional domain serving as a projection of the surface part $D(u + \Delta u) \setminus \Psi(D(u))$. Dissect now the projection of the boundary $D(u + \Delta u) \cap \partial \Omega$ into the small parts $E_i, i = 1, ..., N$ with the maximal diameter not exceeding $\Delta u$. Now taking any point on $E_i$ drow the line orthogonal to the $E_i$, lying on the tangential hyper plane. The set of all such lines set up a surface. We restrict this surface by a such way that the projection of the getting piece was coinside with $\omega'$. It serves as a tiled covering for the surface $D(u + \Delta u) \setminus \Psi(D(u))$.

Let $\bar{\bar{\xi}}_i \in E_i$ be any point, $\rho_i$ be a vector lying on the tangential space to the surface $F = u$ in a point $\bar{\bar{\xi}}_i \in E_i$, orthogonal to $E_i$ with the top in the point $\bar{\bar{\xi}}_i \in E_i$ and with the endpoint at $\bar{\eta}_i$ of the boundary of corresponding piece of tiled covering. At the small $\Delta u$ we have: $|F_i| = |E_i| h_i$ (here $|E_i|$ expresses $n - 2$-dimensional volume of a piece $E_i$), and $h_i = |\rho_i|(1 + o(1))$, i.e. $h_i$ plays a role of height of $F_i$ which approximately we take for the cylindroid with the base $E_i$ (with an error of an order $o(\Delta u)$ for $n - 1$-dimensional volume). Then, we have:

$$\int_{D(u + \Delta u) \setminus D'(u)} \frac{ds}{\|\nabla F\|} = \sum_{i=1}^{N} \int_{(E_i)} \frac{ds}{\|\nabla F\|}(1 + o(1)).$$

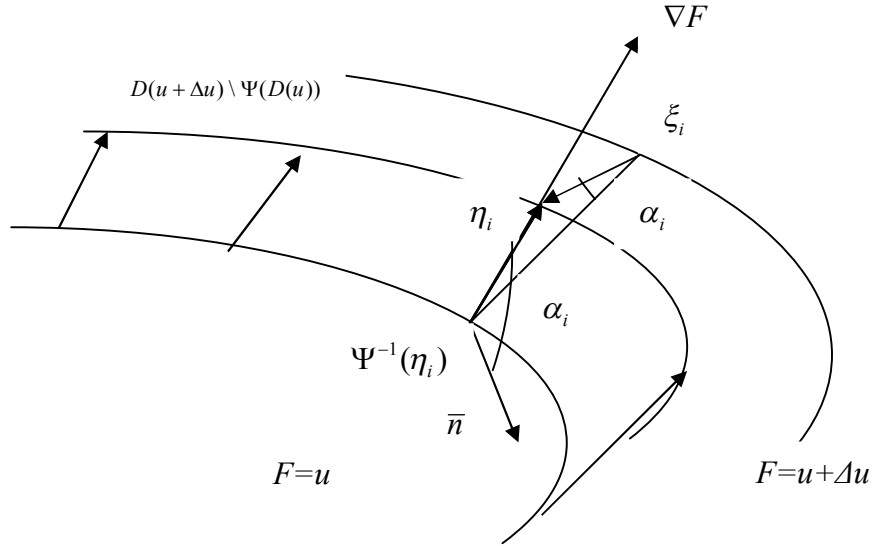

Intersection of tangential hyper planes, respectively, to $\partial \Omega$ and $D(u + \Delta u)$ at the point $\xi_i$, is a tangential $r - 2$- dimensional subspace to $D(u + \Delta u) \cap \partial \Omega$ at the same point. Let's consider three points: a point $\bar{\bar{\xi}}_i$, a point $\bar{\eta}_i$ and a point $\Psi^{-1}(\eta_i)$. Let $\alpha_i$ be an angle between an external normal vector $\bar{n}$ to the boundary of $\Omega$ and a gradient $\nabla F$ (see a drowing). An angle between the segment $[\bar{\eta}_i, \bar{\bar{\xi}}_i]$ and the gradient $\nabla F$, at small $\Delta u$, differs from an angle $\alpha_i$ by a value of $o(1)$ (or their sum is close to $\pi$). From a rectangular triangle we receive (the told above segment $[\Psi^{-1}(\eta_i),\ \bar{\xi}_i]$ is here an hypotenuse):

$$h_i = |\lambda| \cdot \|\nabla F\| ctg\alpha_i(1+o(1)) = \frac{\Delta u}{\|\nabla F\|} ctg\alpha_i(1+o(1)).$$

As $\cos\alpha_i = \bar{n} \cdot \overline{\nabla F} = (\bar{n}, \overline{\nabla F})\ ctg\alpha_i = (\bar{n}, \overline{\nabla F}) / \sqrt{1-(\bar{n}, \overline{\nabla F})^2}$, we have:

$$\int\limits_{D(u+\Delta u)\backslash D'(u)} \frac{ds}{\|\nabla F\|} = \sum_{i=1}^{N} \Delta u \int\limits_{(E_i)} \frac{ctg\alpha_i d\sigma}{\|\nabla F\|^2}(1+o(1)) = \Delta u(1+o(1))\int\limits_{Z} \frac{(\overline{\nabla}, \bar{n})}{\sqrt{1-(\overline{\nabla}, \bar{n})^2}} \frac{d\sigma}{\|\nabla F\|^2},$$

where $d\sigma$ designates $r-2$-dimensional element of volume, and $Z$ denotes an intersection of surfaces $F=u$ and $\partial\Omega$ (it can consist of several pieces). The similar formula is true and for the third surface integral in (4). Therefore, the formula (4) can be rewritten in the view:

$$\varphi'(u) = \lim_{\Delta u \to 0} \frac{\varphi(u+\Delta u)-\varphi(u)}{\Delta u} = -\int\limits_{F(\bar{x})=u} \|\nabla F\|^{-3}(A_1\overline{\nabla}, \overline{\nabla})ds + \int\limits_{Z} \frac{(\overline{\nabla}, \bar{n})}{\sqrt{1-(\overline{\nabla}, \bar{n})^2}} \cdot \frac{ds}{\|\nabla F\|^2}, \quad (5)$$

and the sign before the integral is counted by the scalar product $(\overline{\nabla}, \bar{n})$. The second term on the extreme chain of (5) is possible to represent in the following view using differential forms:

$$\int\limits_{Z} \frac{(\overline{\nabla}, \bar{n})}{\sqrt{1-(\overline{\nabla}, \bar{n})^2}} \cdot \frac{ds}{\|\nabla F\|^2} = \int\limits_{\partial D(\bar{u})} \omega,$$

since for the points which are lying on the boundary and not belonging $\partial\Omega$ we have the relation $(\overline{\nabla}, \bar{n})=0$. Really, use the expression for an area element

$$ds = \sqrt{G} / |J| d\xi,$$

where $\bar{\xi}=(\xi_1,...,\xi_{n-1})$ is a vector of independent variables received after solving of the equation $H(\bar{x})=c$ with regard to one of variables. The polynomial equation $|J|=0$ has a set of solutions of a zero measure in $\Omega$. Therefore, we receive integral

$$\int\limits_{U} \frac{(\overline{\nabla}, \bar{n})}{\sqrt{1-(\overline{\nabla}, \bar{n})^2}} \cdot \frac{\sqrt{G}}{|J|} \frac{d\bar{\xi}}{\|\nabla F\|^2} = \int\limits_{\partial D(\bar{u})} \omega,$$

probably in improper sense (here $U$ denotes the domain of changing for a vector $\bar{\xi}$ of independent variables). We, thus, have received integral of a differential form. Applying Stoks formulae, we will have (see [40, p. 261]):

$$\int_{\partial D(\bar{u})} \omega = \int_{D(u)} d\omega ,$$

where the differential form is defined by the equations defining boundary $\partial\Omega$ of the domin $\Omega$, i.e. by the system of equations of a view $F = u, H = v$. Obviously, from the received integral, that it is possible to pass to the surface integral multiplying by a non-negative function defining an element of the area. From (5) follows:

$$\varphi'(u) = \int_{F(\bar{x})=u} G(\bar{x})ds ,$$

(possible, in improper sense) with some algebraic function $G(\bar{x})$ defined on the domain $\Omega$. Now we will dissect the domain $\Omega$ into union of such subdomains $\Omega_i$, that the function $G(\bar{x})$ keeps its sign invariable. Obviously, at the made assumptions, the function

$$\frac{(\overline{\nabla}, \bar{n})}{\sqrt{1 - (\overline{\nabla}, \bar{n})^2}}$$

is defined for all points of the domain $\Omega$. Then, the integral received above for $\varphi'(u)$ splits into the sum of several surface integrals:

$$\varphi'(u) = \sum \varphi_i'(u), \varphi_i'(\bar{x}) = \int_{\Omega_i, F=u} G(\bar{x})ds , \quad (6)$$

(notice that when we consider the sum of the integrals $\int_{S \subset Z}$ taken on the different sides of the same piece $S$ of a surface, the normal vector $\bar{n}$ changes the sign, and consequently, such a sum is equal to zero); the number of domains on the right part of (6) depends on $\Omega$ and a degree of the polinomial $F$. Let's designate, in a consent with (6) $\varphi(u) = \sum \varphi_i(u) \, \varphi_i(u) = \int_{\Omega_i, F=u} ds / \|\nabla F\|$.

Let's prove the equality $\varphi_i'(u) = \int_{\Omega_i, F=u} G(\bar{x})ds$. The boundary of the domain $\Omega_i$ is defined by the boundary of $\Omega$ and the equation $G(\overline{x}) = 0$. In this case direct application of the reasonings which have been carried out above does not give a desired, but gives a similar representation for a derivative $\varphi'_i(u)$ (we notice that there is much more of such representations). Therefore, the equality we need can be proved by the followimg way.

Let's take the fixed value $u$. On a condition imposed above in each point of boundary of a piece of a surface, lying in $\Omega_i$, the normal vector $\overline{n}$ to the surface $H(\overline{x}) = a$ exists. Let $\delta_i$ designates the boudary of a piece of the surface, concluded in area $\Omega_i$. It is $n-2$-dimensional surface. Let, further, $\overline{x} \in \delta_i$ be any point on the boundary. Then, in this point there is $n-2$-dimensional tangential space and the normal vector $\overline{n}$ considered above to this piece of the surface. Let's take a straight line orthogonal to $n-1$-dimensional subspace generated by tangential space and the normal vector $\overline{n}$. Set of all such straight lines corresponding to every possible points $\overline{x} \in \delta_i$ is a $n-1$-dimensional linear surface. This surface defines (by intersection) the same piece of the surface included in areas $\Omega_i$, and in each point has a normal coinciding with $\overline{n}$. Applying now carried out above reasoning, we receive the required equality $\varphi'_i(u) = \int\limits_{\Omega_i, F=u} G(\overline{x}) ds$.

Thus, the equality $\varphi'(u) = \sum\limits_i \varphi'_i(u) = \sum\limits_i \int\limits_{\Omega_i, F=u} G(\overline{x}) ds$ is true. Since the function under the surface integral does not change its sign, the function $\varphi_i(u)$ is a monotone function. The lemma 5 is proved.

### 4. Proofs of theorems.

*Proof of the theorem 1.* $\mu(\Pi_H)$ is possible to represent in the form of the following surface integral

$$\mu(\Pi_H) = \int\limits_{\Pi \cap E(H)} ds, \qquad (7)$$

where $\Pi$ designates a surface of solutions of the system (3). Let among the parts

$$\Delta_\mu, \mu = 1, \ldots, f, f \le F$$

(see the splitting which has been carried out above) of the surface $\Pi$ the part $\Delta_1$ have a maximal area, and

$$J_1 = \frac{D(f_1, \ldots, f_{n-r})}{D(x_{n-r+1}, \ldots, x_n)}$$

designates the corresponding maximal minor. Let $\Pi_0$ be the area of $\Delta_1$. Then, we receive the following inequality

$$\mu(\Pi_H) \le F\Pi_0.$$

Now we estimate $\Pi_0$.

$$\Pi_0 \le \int_{\Pi_1'} \frac{\sqrt{\sum J_i^2}}{|J_1|} d\xi_1 \ldots d\xi_r \le \binom{n}{r}^{1/2} \int_{\Pi_1'} d\xi_1 \ldots d\xi_r, \quad (8)$$

where $\Pi_1'$ designates the domain of variation of free variables $\xi_1, \xi_2, \ldots, \xi_r$ used for the parametric representation of the surfce $\Pi_1$.

As, the determinant of a matrix $A_0$ is a polynomial, a set where the equality $\det(A_0 \cdot {}^t A_0) = 0$ is satisfied, has a zero Jordan measure. From the expression for the Gram's determinant it follows that

$$|\det D_0|^2 \le \det\left(A_1 \cdot {}^t A_1\right) \le H^2,$$

and $D_0$ is a submatrix of the matrix $A_0$ containing entries Gram's matrix of gradients of which contains the columns of maximal minor $J_1$. Therefore, according to the assumption,

$$|\det J_1|^2 \ge \left(C_{nr}^r\right)^{-1} \det\left(A_1 \cdot {}^t A_1\right) \ge \left(C_{nr}^r\right)^{-1} \left(G_{(1)}\right)^2. \quad (9)$$

Then, for the integral on the right part of (8) we have:

$$\int_{\Pi_1', |\det D_0| \le H} d\xi_1 \cdots d\xi_r \le \sum_{j=1}^{\infty} E_j,$$

where

$$E_j = \int_{2^{-j} H \le |\det D_0| \le 2^{1-j} H} d\xi_1 \cdots d\xi_r.$$

Further reasonings we can carry out as in [47,p.22]. Then we get

$$\int\limits_{\Pi_1',\,|\det D_0|\le H} d\xi_1\cdots d\xi_r \le 2^{r+3}r^{3r}c_0^2\binom{n}{r}^{1/2}\binom{rn}{r}^{1/2}HG_{(1)}^{-1}\wp^r,$$

where

$$\wp = 1 + r^2 + \log\left(rG_{(1)}H^{-r-1}L^{r(r-1)}\right) \le 1 + r^2 + \log r + \log\left\{h\left(G_{(1)}\right)h\left(H\right)^{r+1}h\left(L\right)^{r^2}\right\}.$$

Whent $r \ge 2$ we have

$$r^2\log\left\{h\left(G_{(1)}\right)h\left(H\right)h\left(L\right)\right\} \ge 3r^2 + \log\left\{h\left(G_{(1)}\right)h\left(H\right)^{r+1}h\left(L\right)^{r^2}\right\} \ge \wp.$$

The theorem 1 is proved.

*Proof of the theorem 2.* Let's carry out an induction on $k$. At $k=1$ the theorem 2 follows from the theorem 1. Let's consider passing from the case $k=1$ to the case $k=2$.

At the proof of the theorem 1 we assumed that $G_1(\bar{\bar{\xi}}) > 0$. Differentiation here is carried out with respect to the components of the vector $\bar{\bar{\xi}}$ which we defined from the system of the considered equations and, therefore, lines of the matrix $A_1(\bar{\bar{\xi}})$ are linear combinations of lines of the matrix $A_1(\bar{x})$. Now, we will replace the received estimates by the estimations in term of submatrices of $A_1(\bar{x})$. That gives us to carry out differentiation with respect to independent variables of a vector $\bar{x}$.

The matrix $A_1(\bar{\bar{\xi}})$ is possible to represent as $D(\bar{x})\cdot A_1(\bar{x})$, where $A_1(\bar{x})$ is a matrix entered above, and $D(\bar{x})$ is a matrix of a type:

$$D(\bar{x}) = \begin{pmatrix} 1 & 0 & \cdots & 0 & \varphi_{11} & \cdots & \varphi_{1,4k-r} \\ 0 & 1 & \cdots & 0 & \varphi_{21} & \cdots & \varphi_{2,4k-r} \\ \vdots & \vdots & \ddots & \vdots & \vdots & \ddots & \vdots \\ 0 & 0 & \cdots & 1 & \varphi_{r1} & \cdots & \varphi_{r,4k-r} \end{pmatrix} = \left(E_r \mid \Phi\right),$$

where $E_r$ - an identity matrix of an order $r$, and $\Phi$ is a matrix of the size $r \times (n-r)$. Therefore, any minor, for example, a minor $M_1$ made of the first $r$ columns of a matrix $A_1(\bar{\bar{\xi}})$, is possible to represent in the form of block determinant

$$\delta_1 = \begin{vmatrix} D(\bar{x})\cdot A_1^r(\bar{x}) & 0 \\ \Psi & E_{n-r} \end{vmatrix},$$

and the matrix $A_1^r(\bar{x})$ is a rectangular matrix made of the first $r$ columns of a matrix $A_1(\bar{x})$, $\Psi$ consists of the last $n-r$ lines of a matrix $A_1^r(\bar{x})$. Carrying out elementary transformations over the lines of the last determinant, we present $\delta_1$ in the form

$$\delta_1 = \begin{vmatrix} \left[ A_1^r(\bar{u}) \right] & -\Phi \\ \Psi & E_{n-r} \end{vmatrix}, \quad (10)$$

where the matrix $\left[ A_1^r(\bar{x}) \right]$ is made of the first $r$ lines of a matrix $A_1^r(\bar{x})$ so that two blocks of the first column of determinant $\delta_1$ form a matrix $A_1^r(\bar{x}) : A_1^r(\bar{x}) = \begin{pmatrix} \left[ A_1^r(\bar{x}) \right] \\ \Psi \end{pmatrix}$.

Let's dissect the domain $\Pi_1'$ into two parts: in the first the condition $\sqrt{G_1(\bar{\bar{\xi}})} \geq G_{(1)}$ is satisfied, for some positive $G_1$, and in the remaining part of the domain $\Pi_1'$ the inequality $\sqrt{G_1(\bar{\bar{\xi}})} \leq G_{(1)}$ is satisfied. Designating by $\mu_1$ and $\mu_2$ the areas of the corresponding parts of the surface $\Pi_1$, for our surface integral we receive:

$$\int_{\Pi_1, |\det D_0| \leq H} d\sigma \leq \mu_1 + \mu_2 . \ (11)$$

For the estimation of $\mu_1$ we use the relation (9). We have:

$$\mu_1 \leq 2^{r+3} r^{3r} c_0^2 F T_0 \binom{n}{r}^{3/2} \binom{nr}{r}^{1/2} H G_{(1)}^{-1} \wp_1^r ; \ (12)$$

$$\wp_1 = r^2 \log \left\{ h\left( G_{(1)} \right) h(H) h(L) \right\} .$$

The estimation of $\mu_2$ can be reduced to the estimation similar to (12). For do it at first we allocate area where the condition of a view $\eta < \sqrt{G_1(\bar{\bar{\xi}})} \leq 2\eta \leq G_{(1)}$ (the corresponding we designate as $\mu_2^{(1)}$) is satisfied. Let $M_1 = M_1(\bar{\bar{\xi}})$ be a minor of a matrix $A_1(\bar{\bar{\xi}})$ containing elements, gradients of which compile the maximal minor of a matrix $A_2(\bar{\bar{\xi}})$. It obviously, that $|M_1| \leq \sqrt{G_1(\bar{\bar{\xi}})} \leq 2\eta$. Further, replacing the condition $\eta < \sqrt{G_1(\bar{\bar{\xi}})} \leq 2\eta$ by a condition $|M_1| \leq 2\eta$, we will increase only the considered integral. Therefore,

$$(2\eta)^{-1} \mu_2^{(1)} \leq \int_{\Pi_1', |M_1| \leq 2\eta} \frac{1}{\sqrt{G_1\left( \bar{\bar{\xi}} \right)}} d\bar{\xi} = c_0 \int_{\Pi_1', |M_1| \leq 2\eta} d\bar{\xi} \int_{\| D_1 \bar{v} \| \leq 1} d\bar{v} ,$$

thus the integral on the right part undertaken over the part of the product $\Pi_1' \times \mathbf{R}^n$

where the imposed conditions on variables are satisfied. At every point $\bar{x} \in \Pi_1$, all functions depending on $\bar{x}$ contiued, as constants, by parallel translation by vectors from $\mathbf{R}^n$ orthogonal to the surface $\Pi_1$. Let $\pi'_1(\bar{v})$ is a part of the surface $\Pi_1$ for the points of which, at every $\bar{v}$ conditions $|M_1| \le 2\eta$ and $\|D_1\bar{v}\| \le 1$ are satisfied. Having changed the orders of integration, we receive:

$$(2\eta)^{-1}\mu_2^{(1)} \le 2c_0 \int\limits_{\|D_1\bar{v}\| \le 1} d\bar{v} \int\limits_{\pi'_1(\bar{v}), \eta \le |M_1| \le 2\eta} d\sigma . \quad (13)$$

Let's make in the internal surface integral the exchange of variables by the formulas: $\bar{\beta} = D_1\bar{v}$ applying the lemma 1, [47]. Then in designations of this lemma:

$$\int\limits_{\pi'_1(\bar{v}), \eta \le |M_1| \le 2\eta} \sqrt{G}\, \frac{d\sigma}{\sqrt{G}} = \int\limits_{\omega(\eta)} |\det Q| \sqrt{G}\, \frac{d\sigma'}{\sqrt{G'}} , \quad (14)$$

and $\omega(\eta)$ is a prototype of the surface at specified mapping, and

$$G' = \det(JQ \cdot {}^t Q^t J) ;$$

further, $Q$ is a Jacoby matrix of change of vrables which is equal to the inverse for a matrix

$$\frac{\partial \bar{\beta}}{\partial \bar{x}} = \frac{\partial(\beta_1, ..., \beta_n)}{\partial(x_1, ..., x_n)} ,$$

and $J = A_0$, $0.5 \le \sqrt{G} \le 1$ (in agree with (3)). It is known that the lines of the matrix $A_0$ form a subspace $M$ orthogonal to the subspace $M'$ being span of the lines of the matrix $D(\bar{x})$. Therefore, $R^n = M \oplus M'$, and it is possible to present an element of volume as $d\bar{w} = d\bar{y}d\bar{z}$. Each vector $\bar{w} \in R^n$ uniquely written in the form of the sum of vectors $\bar{y}$ and $\bar{z}$ from the subspaces $M$ and $M'$, thus $\|Q(\bar{y}+\bar{z})\| \le \|Q\bar{y}\| + \|Q\bar{z}\|$. Let $\bar{y}_1, ..., \bar{y}_r$ be the basis consisting of lines of a matrix $A_0$, $\bar{z}_1, ..., \bar{z}_{n-r}$ is a basis consisting of lines of a matrix $D = D(\bar{x})$. Making change of variables in the integral given below by formulas $\bar{w} = W\bar{x}$ where $W = \begin{pmatrix} A_0 \\ D \end{pmatrix}$, we receive

$$|\det Q^{-1}| = c_0 \int\limits_{\|Q\bar{w}\| \le 1} d\bar{w} = c_0 \int\limits_{\|Q(A_0\bar{u}+D\bar{v})\| \le 1} \frac{1}{\sqrt{G}\sqrt{D \cdot {}^t D}} d\bar{u}d\bar{v} \ge$$

$$\geq c_0 \int\limits_{\|Q A_0 \bar{u}\| \leq 1/2} \frac{1}{\sqrt{G}} d\bar{u} \int\limits_{\|QD\bar{v}\| \leq 1/2} \frac{1}{\sqrt{D \cdot {}^t D}} d\bar{v} \, .$$

Therefore,

$$\left| \det Q^{-1} \right| = \frac{\Gamma(1+n/2)}{\pi^{n/2}} \int\limits_{\|Q\bar{w}\| \leq 1} d\bar{w} \geq$$

$$\geq \frac{2^{-n} \Gamma(1+n/2)}{\Gamma(1+r/2) \Gamma(1+(n-r)/2) \sqrt{D \cdot {}^t D} \sqrt{G}} \frac{1}{\sqrt{G'}} \frac{1}{\sqrt{|{}^t D' Q Q D|}} \, .$$

Substituting in (13), we find:

$$\int\limits_{\pi'_1(\bar{v}), \eta \leq |M_1| \leq 2\eta} d\sigma \leq c' \int\limits_{\omega(\eta)} \sqrt{D \cdot {}^t D} \sqrt{{}^t D' Q Q D} \cdot G d\sigma' \, ,$$

from some positive constant $c'$. Noting that every enties of the matrix $\Phi$, has a modulus does not exceeding 1, and that $\sqrt{D \cdot {}^t D} \leq (n-r)^r \sqrt{G} \leq 1$, we receive:

$$\int\limits_{\pi'_1(\bar{v}), \eta \leq |M_1| \leq 2\eta} d\sigma \leq (n-r)^r c' \int\limits_{\omega(\eta)} \sqrt{{}^t D' Q Q D} \cdot d\sigma' \, .$$

We have:

$$\mu_2^{(1)} \leq 4(n-r)^r c' \eta \int\limits_{\|D_1\bar{v}\| \leq 1} d\bar{v} \int\limits_{\omega(\eta)} \sqrt{{}^t D' Q Q D} \cdot d\sigma' \, .$$

Let's consider restriction of biyectiv mapping $\bar{\beta} = D_1 \bar{v}$ to the variety of solution of the system. The matrix $QD$, in each point $\bar{x} \in \Pi'_1$ being a solution of the system, is a Jacoby matrix of invers transformation, i.e. its invers coincides with a Jacoby matrix of the change of varables $\bar{\beta} = D_1 \bar{v}$ *taken with reqard to independent variables of* $\bar{\bar{\xi}}$. Then, it is possible return to (8), and further, using the reasonings have been carried out above, we come to the relation of the type already received above in (12), with replacement $\delta_1$ by $\delta_2$ and $H$ by $\eta$ (note that gradients of elements of the matrix $M_1$ contain columns of the maximal minor of a matrix $A_2(\bar{\bar{\xi}})$). So,

$$\mu_1 \leq T_1 H G_{(1)}^{-1} \wp_1^r \, ;$$

$$T_0 F \cdot 2^{r+3} r^{3r} c_0^2 \binom{n}{r}^{3/2} \binom{nr}{r}^{1/2} \le T_1 = 2^{r+3} n^{2r} r^{7r/2} c_0^2 T_0 F \,,$$

$$\mu_2^{(1)} \le 4(n-r)^r c' T_1 \eta G_{(2)}^{-1} \wp_2^r \,; \quad \wp_1 = r^2 \log \left\{ h\left(G_{(1)}\right) h(H) h(L) \right\},$$

$$\wp_2 = r^2 \log \left\{ h\left(G_{(1)}\right) h\left(G_{(2)}\right) h(L) \right\}, \,.$$

Substituting consequently $\eta/2, \eta/4,\dots$ for $\eta = G_{(1)}$, and summarising the received estimates, we find:

$$\mu_2 \le 4(n-r)^r c'' T_1 \cdot G_{(1)} G_{(2)}^{-1} \wp_2^r \le 2^{r+5} n^{3r} r^{7r/2} c'' T_0 F \cdot G_{(1)} G_{(2)}^{-1} \wp_2^r, \,(15)$$

where $c''$ is a sum of series like considered above. It, therefore, depends only on $n, r$ and the degree of the polynomial. Writing $G_{(1)} = H^{1/2} G_{(2)}^{1/2}$, from (11), (12) and (15) we find the estimation:

$$\int\limits_{\Pi'_1, |\det D_0| \le H} d\xi_1 \cdots d\xi_r \le \left(1 + 4n^r c''\right) c_0^2 2^{r+3} \left(nr^2\right)^{2r} T T_0 H^{1/2} G_{(2)}^{-1/2} \left(\wp''\right)^r \,;$$

$$\wp'' = 3r^2 \log H_0; H_0 = \max \left\{ h(H), h\left(G_{(1)}\right), h\left(G_{(2)}\right), h(L) \right\}.$$

Thus, we have passed from the estimation for the case $k=1$ to the case $k=2$. Repeating similar reasonings, with insignificant changes, we can using an estimation for the case $k$ to pass to the case $k+1$. The theorem 2 is proved.

*Proof of the theorem 3.* For the proof of the theorem 3 it is enough to prove that for any natural $k$ we have an inequality $G_{(k)} \ge G_k$. In the entered designations $G_{(k)} = \min \det \left( A_{k\bar{\xi}} \cdot {}^t A_{k\bar{\xi}} \right)$. For any point $\bar{x} = \bar{x}\left(\bar{\bar{\xi}}\right)$ lying on the variety of solutions of this system

$$\det \left( A_{k\bar{\xi}} \cdot {}^t A_{k\bar{\xi}} \right)^{-1/2} = c_0 \int\limits_{\left\| {}^t A_{k\bar{\xi}} \bar{u} \right\| \le 1} d\bar{u} = c_0' \int\limits_{\bar{x} \in \pi, \left\| {}^t A_k \bar{x} \right\| \le 1} ds \,,$$

where the surface integral taken on the tangential space to the surface at the point $\bar{x} = \bar{x}\left(\bar{\bar{\xi}}\right)$. Considering tangential space as a subspace $R^n$ of a dimension $r$, and taking (see [13, p. 148]) the maximal value of the surface integrals over all $r$-dimensional subspaces, we will receive the necessary estimation. The proof of the theorem 3 is complete.

*Proof of the theorem 4.* Let's present this integral in the view

$$J = J_1 + J_2,$$

where

$$J_1 = \int_{\Omega_1} e^{2\pi i F(\overline{x})} d\overline{x}, \quad J_2 = \int_{\Omega_2} e^{2\pi i F(\overline{x})} d\overline{x},$$

and $\Omega_1$ is defined by a condition $\|\nabla F\| \le H$, at some $H$ defined below. The domain $\Omega_1$ consists of finite number of subdomains. The same it is possible to state about the area $\Omega_2$ which is defined by an inequality $\|\nabla F(\overline{x})\| \ge H$.

To estimate $J_1$ we will make at first the chage of variables defined as below:

$$u_1 = \frac{\partial F}{\partial x_1}, ..., u_r = \frac{\partial F}{\partial x_r}.$$

Then we receive:

$$\left| \int_{\Omega_1} e^{2\pi i F(\overline{x})} dx \right| \le \int_{\substack{\sqrt{u_1^2 + ... + u_r^2} \le H \\ \overline{x} \in \Omega}} \frac{1}{\left| \det A_1(\overline{x}) \right|} du_1 ... du_r,$$

where $\overline{x}$ is considered as a function of $\overline{u} = (u_1, ..., u_r)$. The part of this integral defined by the condition $\left| \det A_1(\overline{x}) \right| \le V$ we estimate by using of the theorem 2 as a value (taking $r = n$ in the theorem 2):

$$\le CV^{\frac{1}{k-2}} G^{-\frac{1}{k-2}} \wp^r,$$

where $C$ depends only on $degF$, $r$, $k$, and $G^2$ is a minimal value of the determinant of a matrix $A_{k-1}{}' A_{k-1}$. Remaining part of the integral is estimated as a value:

$$\le V^{-1} \int_{\sqrt{u_1^2 + ... + u_r^2} \le H} du_1 ... du_r \le c_1'' H^r V^{-1},$$

and $c_1''$ depends only on $r$. Now writing

$$V = H^{\frac{r(k-2)}{k-1}} G^{\frac{1}{k-1}}$$

we have:

$$\left| \int_{\Omega_1} e^{2\pi i F(\overline{x})} dx \right| \le C_1 H^{\frac{r}{k-1}} G^{-\frac{1}{k-1}} \wp_k^{r+1}, \quad (16)$$

and $C_1$, $r$, and $k$ depends only on $degF$.

Let's transform now integral $J_2$. From the consequence to the lemma 1 of [46] we have:

$$\int_{\Omega_2} e^{2\pi i F(\bar{x})} d\bar{x} = \int_m^M e^{2\pi i u} du \int_{F(\bar{x})=u} \frac{ds}{\|\nabla F\|} = \int_m^M \left( \int_{F(\bar{x})=u} \frac{ds}{\|\nabla F\|} \right) (\cos 2\pi u + i \sin 2\pi u) \, du .$$

Here $m$ and $M$, respectively, the minimal and maximl values of the function $F(\bar{x})$ which are reached in $\Omega$. The internal integral is the surface integral taken over the surface $F(\bar{x}) = u$. Integral

$$\varphi(u) = \int_{F(\bar{x})=u} \frac{ds}{\|\nabla F\|}$$

as the function depending on $u$, according to Lemma 1, is possible to present in the form of the sum of no more than $K_0$ monoton functions. Then for $J_2$ we receive an estimation:

$$|J_2| \leq 2K_0 \int_{\Pi} \frac{ds}{\|\nabla F\|}, \quad (17)$$

and $\Pi$ means that piece of the surface $F(\bar{x}) = u$ for which the surface integral accepts the maximal value. Let the maximum on the right part of (16) is reached in the point $u = u_0$ (for simplicity the area of this piece also we designate as $\Pi$). Let's dissect the domain $\Omega$ into finite number of parts where one of partial derivatives $\partial F / \partial x_i$ accepts greatest values for its modulus among all partial derivatives of the first order. Let, for example $i = 1$. Then, replacing the condition $U \leq \|\nabla F\| \leq 2U$ in the surface integral (17) by the condition

$$\sqrt{\left( \partial F / \partial x_2 \right)^2 + \cdots + \left( \partial F / \partial x_r \right)^2} \leq 2U ,$$

we only will increase this integral. We have:

$$\int_{U \leq \|\nabla F\| \leq 2U} \frac{ds}{\|\nabla F\|} \leq U^{-1} \int_{\sqrt{(\partial F / \partial x_2)^2 + \cdots + (\partial F / \partial x_r)^2} \leq 2U} ds \leq \Pi U^{-1} , \quad (18)$$

where $\Pi$ is a maximal value of the area of a considered surface $F(\bar{x}) = u_0$. As

$$ds = \frac{\|\nabla F\|}{|\partial F / \partial x_1|} dx_2 \cdots dx_r \leq \sqrt{r} \, dx_2 \cdots dx_r ,$$

then for the same surface integral is possible to receive another estimate making the change of variables of a view $u_i = \partial F / \partial x_i; i = 2, ..., r$ :

$$\int\limits_{U \leq \|\nabla F\| \leq 2U} \frac{ds}{\|\nabla F\|} \leq \sqrt{r} U^{-1} \int\limits_{\sqrt{u_2^2 + \cdots + u_r^2} \leq 2U} \Delta^{-1} dx_2 \cdots dx_r \leq c_2 U^{r-2} G_2^{-1} . \ (19)$$

Here we used that the determinant

$$\Delta = \left| \det \begin{pmatrix} \partial^2 F / \partial x_2^2 & \cdots & \partial^2 F / \partial x_2 \partial x_r \\ \vdots & \ddots & \vdots \\ \partial^2 F / \partial x_r \partial x_2 & \cdots & \partial^2 F / \partial x_r^2 \end{pmatrix} \right|$$

it is not less than the product of $r-1$ smallest singulyar numbers of the following matrix of an order $r : \left( \partial^2 F / \partial x_i \partial x_j \right)_{1 \leq i, j \leq r}$. Let's dissect a piece of the surface, lying in the $\Omega_2$ into two parts so that in one of them $\Delta \leq R$, and in another $\Delta > R$, for some positive $R$ the exact value of which we will pick up later. For the surface integral taken on the second piece, we have:

$$\int\limits_{\Delta > R} \frac{ds}{\|\nabla F\|} \leq \sqrt{r} U^{-1} \int\limits_{\sqrt{u_2^2 + \cdots + u_r^2} \leq 2U} dx_2 \cdots dx_r \leq c_2 U^{r-2} R^{-1} . \ (20)$$

Integral on the first piece we estimate by using of the theorem 3, taking $A_0$ equal to the determinant of the matrix $\Delta$ :

$$\int\limits_{\Delta \leq R} \frac{ds}{\|\nabla F\|} \leq c_3 \sqrt{r} U^{-1} R^{1/(k-2)} G_{k-2}^{-1/(k-2)} \wp_k^r .$$

Let's determine $R$ by a condition: $U^{-1} R^{1/(k-2)} G_{k-2}^{-1/(k-2)} = U^{r-2} R^{-1}$ i.e $R = U^{\frac{(r-1)(k-2)}{k-1}} G_{k-2}^{1/(k-1)}$. Then we receive an estimation

$$|J_2| \leq c_4 U^{\frac{r-k}{k-1}} G_{k-2}^{-\frac{1}{k-1}} \wp_k^r . \ (21)$$

Then, noting that the segment $[H_1, H_0]$ can be dissected into no more, than $1 + \left[ \log_2 \left( H_0 H_1^{-1} \right) \right]$ segments of a view $[U, U'], U' \leq 2U$ , we receive for $k \leq r$

$$|J_2| \leq c_4 H_0^{\frac{r-k}{k-1}} G_{k-2}^{-\frac{1}{k-1}} \wp_k^{r+1} . \ (22)$$

If $k > r$ , counting the bound (16) at $H = (GG_{k-1}^{-1})^{\frac{1}{k}}$ where $G^2$ is the minimal value of the determinant of the matrix $A_{k-1} \cdot {}^t A_{k-1}$ , finally, for $J$ we receive the follow-

ing estimation

$$J \leq C(r,k,F)K_0 \cdot G^{\frac{r-k}{k(k-1)}} G_{k-1}^{-\frac{r}{k(k-1)}} \wp_k^{r+1};$$

here $C = C(r,k,F)$ depends only on $r$, $k$ and $F$.

If now $\Omega$ is a cube with an edge 1, as it is clear from the following inequality now in (18) $\Pi$ is possible to replace by 1:

$$\int_{\sqrt{u_2^2 + \cdots + u_r^2} \leq 2U} dx_2 \cdots dx_r \leq \int_0^1 \cdots \int_0^1 dx_2 \cdots dx_r = 1.$$

Then, applying (22) for $\bar{x}$ such that $\|\nabla F\| \leq U$ we receive taking $U = G_{k-2}$:

$$|J| \leq c_4 G_k^{-\frac{1}{r-1}} \wp_k^{r+1}.$$

It proves 2), b) of the theorem 4. The theorem 4 is proved.

*Consequence.* The estimation below takes place:

$$\left| \int_0^1 \cdots \int_0^1 e^{2\pi i F(\bar{x})} d\bar{x} \right| \leq c_4 K \lambda^{-1} \wp_k^{r+1},$$

where $c_4 = c\,(r,\,k,\,F)$ and $K$ are positive constants, and $\lambda$ is defined as below:

1) for $k \leq r$ we put $\lambda = \min_{\bar{x} \in \Omega} (\lambda_1 + \ldots + \lambda_r)$, where $\lambda_i$ is a minimal singular number of the matrix $A_{i-1}$, $1 \leq i \leq r$;

2) for $k > r$, we put $\lambda = \min_{\bar{x} \in \Omega} \lambda_k^{r/k}$, where $\lambda_k$ is a minimal singular number of the matrix $A_{i-1}$.

The first statement is received in [45]. Both results are easy consequences of the theorem 4.

Haydar Aliev avenue, 187,

Ganja Stat University,

Ganja, Azerbaijan.